\documentclass[11pt,reqno]{amsart}

\usepackage[T1]{fontenc}
\usepackage{lmodern}
\usepackage{amsmath,amssymb,amsthm,mathtools}
\usepackage{booktabs}
\usepackage[hidelinks]{hyperref}

\usepackage{graphicx}
\usepackage{tikz}
\usetikzlibrary{arrows.meta,decorations.markings}

\theoremstyle{plain}
\newtheorem{theorem}{Theorem}[section]
\newtheorem{proposition}[theorem]{Proposition}
\newtheorem{lemma}[theorem]{Lemma}
\newtheorem{corollary}[theorem]{Corollary}

\theoremstyle{definition}
\newtheorem{definition}[theorem]{Definition}

\theoremstyle{remark}
\newtheorem{remark}[theorem]{Remark}

\numberwithin{equation}{section}

\newcommand{\C}{\mathbb C}
\newcommand{\R}{\mathbb R}
\newcommand{\Z}{\mathbb Z}
\newcommand{\PP}{\mathbb P}
\newcommand{\OO}{\mathcal O}
\newcommand{\Amo}{\mathcal A}
\newcommand{\ind}{\mathbf 1}
\DeclareMathOperator{\CL}{CL}
\DeclareMathOperator{\Res}{Res}
\DeclareMathOperator{\Log}{Log}
\DeclareMathOperator{\Int}{int}
\DeclareMathOperator{\Conv}{conv}
\DeclareMathOperator{\vol}{vol}

\title[Tropical limits of holomorphic linking]
{Tropical limits of holomorphic linking for toric hypersurfaces and points}

\author{Yuguang Zhang}
\address{Xi'an Jiaotong--Liverpool University, Suzhou, China}
\email{Yuguang.Zhang@xjtlu.edu.cn}
\thanks{Supported in part by NSFC grant 12531001.}


\begin{document}

\begin{abstract}
A holomorphic top form on a smooth hypersurface in a smooth projective
toric variety has a unique meromorphic residue lift when the ambient
dimension is at least two. Evaluating this lift at a point, relative
to the logarithmic volume form, gives a complex linking number.
For a one-parameter family of hypersurfaces and a normalised monomial
residue, we prove that this function converges to the indicator of
the corresponding tropical chamber. The convergence is uniform over
compact sets disjoint from the tropical hypersurface, independently
of the arguments of the point coordinates. When the chamber is
nonempty and its index is an interior lattice point, the limit is the
classical linking number of the oriented chamber boundary with the
limiting logarithmic point. A plane cubic illustrates both this
interpretation and its failure on a tropical wall.
\end{abstract}

\maketitle

\section{Introduction}

Let \(\Gamma\) be the positively oriented boundary of a polygon in
the plane. Its winding number about a point is \(1\) inside the
polygon and \(0\) outside. There is no comparable notion of
``inside'' for a complex curve in a complex surface: the curve has
real codimension two, not one. Nevertheless, a curve equipped with
a holomorphic differential and a point outside it can determine a
complex number. The question in this note is how such a number
can recover an inside--outside distinction.

The simplest example already exhibits the answer. Consider the cubic
\[
  1+z_1^3+z_2^3+t z_1z_2=0,\qquad t>1.
\]
For a natural choice of differential on its projective closure, the
linking  at a point \(z\) outside the curve is
\[
  \frac{t z_1z_2}{1+z_1^3+z_2^3+t z_1z_2}.
\]
Write \(x_i=\log|z_i|/\log t\). If
\[
  1+x_1+x_2>\max\{0,3x_1,3x_2\},
\]
the numerator dominates the other terms in the denominator and the
quotient tends to \(1\). These inequalities describe the interior
of a triangle in the \(x\)-plane. In each region where one of the
other three terms strictly dominates, the quotient instead tends
to \(0\). Thus it is the boundary of a \emph{real triangle in
logarithmic coordinates}, not the complex cubic itself, that
encloses the point.

Theorem~\ref{thm:main} establishes this interpretation for toric
hypersurfaces in arbitrary dimension. Two ingredients have different
roles. The residue construction explains why the monomial quotient
is a linking function. The tropical chamber explains why its limit
is a classical linking number. Once these objects have been
identified, a uniform estimate follows by comparing the largest
monomial with the second largest.

The construction is motivated by polar homology and polar linking
as developed by Khesin and Rosly \cite{KR,KRlink}. For holomorphic
linking of curves in Calabi--Yau threefolds, see Frenkel and Todorov
\cite{FT}. See \cite{FT,KR,KRlink} and the references therein for the project exploring the analogues between topology and complex geometry.    On the tropical
side, the relevant background is the relation between hypersurfaces,
amoebas, and polyhedral decompositions \cite{Mik}. Leung and Wan
\cite{LW} study the behaviour of holomorphic forms in degenerating
hypersurfaces; our concern is the evaluation of their meromorphic
lifts at points \emph{off} the hypersurface. 

\section{Linking by degree and by residue}

\subsection{The real linking number}

Let \(\Xi\subset\R^n\) be a compact domain with piecewise smooth
boundary \(\Gamma=\partial\Xi\), oriented by the outward-boundary
convention. For \(q\notin\Gamma\), define
\[
  L(\Gamma,q)
  =\deg\left(\Gamma\longrightarrow S^{n-1},\quad
             x\longmapsto\frac{x-q}{|x-q|}\right).
\]
Equivalently, pulling back the normalised volume form of the sphere
gives
\[
  L(\Gamma,q)=
  \frac{1}{\vol(S^{n-1})}
  \int_\Gamma
 \sum_{i=1}^{n} \frac{\displaystyle(-1)^{i-1}(x_i-q_i)
       \,dx_1\wedge\cdots\wedge\widehat{dx_i}\wedge\cdots\wedge dx_n}
       {|x-q|^n}.
\]
For \(n=2\) this is the winding number. In any dimension,
\begin{equation}\label{eq:classical}
  L(\Gamma,q)=
  \begin{cases}
    1,&q\in\Int(\Xi),\\
    0,&q\notin\Xi.
  \end{cases}
\end{equation}
Indeed, if \(q\notin\Xi\), the radial map extends over \(\Xi\).
If \(q\in\Int(\Xi)\), remove a small open ball centred at \(q\).
The radial map extends over what remains, so its degree on
\(\Gamma\) is its degree on the positively oriented small
sphere, namely \(1\).

\subsection{The holomorphic linking number}

Let \(Z\) be a compact complex manifold of dimension \(n\), let
\(i:X\hookrightarrow Z\) be a smooth hypersurface, and let \(K_Z\)
and \(K_X\) denote the canonical bundles. A section of \(K_Z(X)\)
is a meromorphic top form on \(Z\) with at most a simple pole
along \(X\) and no poles elsewhere. The residue sequence is
\[
  0\longrightarrow K_Z\longrightarrow K_Z(X)
   \xrightarrow{\Res_X}i_*K_X\longrightarrow0.
\]
Our convention is
\[
  \Res_X\left(\frac{ds}{s}\wedge\alpha+\beta\right)=\alpha|_X,
\]
where \(s=0\) is a local equation for \(X\) and \(\alpha,\beta\)
are holomorphic.

Suppose that
\begin{equation}\label{eq:vanishing}
  H^0(Z,K_Z)=H^1(Z,K_Z)=0.
\end{equation}
The long exact cohomology sequence then identifies
\(H^0(Z,K_Z(X))\) with \(H^0(X,K_X)\). In particular, each
holomorphic top form \(\eta\) on \(X\) has a unique
\emph{residue lift}: a form \(\Theta\) on \(Z\) with
\(\Res_X\Theta=\eta\).     $(X, \eta) $ is the polar boundary of $(Z, \Theta)$  in the sense of the polar homology of  \cite{KR,KRlink}. 

\begin{definition}\label{def:linking}
Assume \eqref{eq:vanishing}, and let \(\Omega\) be a meromorphic
top form on \(Z\) with no zeros.   For \(p\in Z\setminus X\) and
\(\tau\in\C\), define the holomorphic linking number by
\[
  \CL_\Omega((X,\eta),(p,\tau))
  =\tau\,(\Theta\,\Omega^{-1})(p),
\]
where \(\Theta\) is the residue lift of \(\eta\).
\end{definition}

At a point where \(\Omega\) is regular this is just
\(\tau\Theta(p)/\Omega(p)\), the ratio of two elements of the same
one-dimensional vector space. No separate definition is needed at
a pole of \(\Omega\): its inverse is a holomorphic section of
\(K_Z^{-1}\) and vanishes there. The choice of \(\Omega\) matters;
rescaling it rescales the linking function inversely. We take
\(\tau=1\) below.

The case of  $n=2$ has been used in the definition of the holomorphic linking number of two curves in Calabi-Yau threefolds in \cite{FT,KR}. 

\section{Monomial residues on a toric hypersurface}

Let \(M=\Z^n\), with \(n\ge2\), and let \(A\subset M\) be a finite
set whose convex hull \(\Delta\) has dimension \(n\). We do not
require \(A\) to contain every lattice point of \(\Delta\). Assume
that \(\Delta\) has a smooth normal fan. Equivalently, at every vertex
the primitive vectors directed along the edges into \(\Delta\)
form a basis of \(M\). Let \(Z_\Delta\) be the associated smooth
projective toric variety and \(T=(\C^*)^n\) its dense torus. We
use the conventions of \cite{CLS}.

For \(k=(k_1,\ldots,k_n)\in M\), write
\(z^k=z_1^{k_1}\cdots z_n^{k_n}\). On \(T\) set
\[
  \Omega=\frac{dz_1}{z_1}\wedge\cdots\wedge\frac{dz_n}{z_n}.
\]
It extends to \(Z_\Delta\) with a simple pole along each irreducible
component of the toric boundary and no zeros. Moreover,
\(H^j(Z_\Delta,\OO_{Z_\Delta})=0\) for \(j>0\) \cite{CLS}.
By Serre duality, \(H^0(Z_\Delta,K_{Z_\Delta})\) and
\(H^1(Z_\Delta,K_{Z_\Delta})\) therefore vanish. Here the
restriction \(n\ge2\) ensures that \(n-1>0\).

Let
\[
  f(z)=\sum_{k\in A}c_kz^k,\qquad c_k\in\C^*,
  \qquad A^\circ=A\cap\Int(\Delta).
\]
The polynomial defines a section of the line bundle associated to
\(\Delta\). Since its support contains all vertices of \(\Delta\),
its zero divisor has no toric boundary component. Denote this
divisor by \(X\), and assume it is smooth.

\begin{proposition}[Monomial residue lifts]\label{prop:lift}
For \(m\in A^\circ\), put
\[
  \Theta_m=\frac{c_mz^m}{f}\,\Omega,
  \qquad
  \eta_m=\Res_X\Theta_m.
\]
Then \(\Theta_m\) has at most a simple pole on \(X\) and no other
poles. In particular, \(\eta_m\) is holomorphic on \(X\), and
\(\Theta_m\) is its unique residue lift. For \(p\in T\setminus X\),
\begin{equation}\label{eq:quotient}
  \CL_\Omega((X,\eta_m),(p,1))
  =\frac{c_m p^m}{\sum\limits_{k\in A}c_kp^k}.
\end{equation}
\end{proposition}

\begin{proof}
Only the behaviour at the toric boundary needs to be checked. Fix a vertex
\(\nu\) of \(\Delta\), and let \(u_1,\ldots,u_n\) be the primitive
edge vectors directed from \(\nu\) into \(\Delta\). The monomials
\(w_i=z^{u_i}\) are coordinates on the vertex chart
\(U_\nu\simeq\C^n\); the toric boundary in this chart is
\(w_1\cdots w_n=0\).

Write each exponent \emph{relative to the vertex}:
\[
  k-\nu=\sum_{i=1}^n d_i(k)u_i .
\]
The tangent cone at \(\nu\) contains \(\Delta-\nu\), so
\(d_i(k)\ge0\). These coordinates are integers because the \(u_i\)
are a lattice basis. If \(k=m\) is interior to \(\Delta\), then
\(d_i(m)\ge1\) for every \(i\): equality to zero would put \(m\)
on a facet through \(\nu\).

The regular local equation for \(X\) in this chart is
\[
  F_\nu(w)=z^{-\nu}f(z)
  =\sum_{k\in A}c_k\prod_{i=1}^n w_i^{d_i(k)}.
\]
The change of lattice basis has determinant \(\pm1\). Consequently,
\[
  \Theta_m
  =\pm c_m\,
    \frac{\prod\limits_{i=1}^n w_i^{d_i(m)-1}}{F_\nu(w)}
    \,dw_1\wedge\cdots\wedge dw_n .
\]
Every exponent in the numerator is nonnegative. This is precisely
where the interior-point assumption is used: it cancels the poles
of the logarithmic volume form at the toric boundary. The vertex
charts cover \(Z_\Delta\), so \(\Theta_m\) is a section of
\(K_{Z_\Delta}(X)\). The residue sequence gives the assertions about
\(\eta_m\) and uniqueness. Evaluating on \(T\) gives
\eqref{eq:quotient}.
\end{proof}

The coefficient \(c_m\) in the numerator fixes the normalisation:
the quotient is the \(m\)-th term of \(f\) divided by the entire
sum. In particular, multiplying all coefficients by a common
nonzero scalar changes neither \(\Theta_m\) nor its evaluation.

\section{Tropical chambers and the linking limit}

\subsection{The logarithmic coordinates of a family}

Fix coefficients \(a_k\in\C^*\) and a height function \(v:A\to\R\).
For \(t>1\), set
\begin{equation}\label{eq:family}
  f_t(z)=\sum_{k\in A}a_kt^{-v(k)}z^k,\qquad
  \Log_t(z)=
  \frac{1}{\log t}\bigl(\log|z_1|,\ldots,\log|z_n|\bigr).
\end{equation}
Let \(X_t\) be the zero divisor of \(f_t\) in \(Z_\Delta\).
Throughout we assume that \(X_t\) is smooth for all sufficiently
large \(t\).  For \(x=\Log_t(z)\), the absolute value of the \(k\)-th term of
\(f_t(z)\) is
\begin{equation}\label{eq:size}
  |a_kt^{-v(k)}z^k|=|a_k|\,t^{\ell_k(x)},
  \qquad \ell_k(x)=\langle k,x\rangle-v(k).
\end{equation}
Thus the largest exponent, rather than the largest coefficient,
determines the dominant monomial as \(t\) grows.

\begin{figure}[htbp]
\centering
\begin{tikzpicture}[
    x=0.68cm, y=0.68cm,
    font=\footnotesize,
    line cap=round, line join=round
]

\begin{scope}
  \filldraw[fill=teal!55, draw=teal!70!black, thick]
    (0,0)
    .. controls (1.333,0) and (2.667,0.667) .. (4,2)
    .. controls (4.067,1.867) and (4.067,1.767) .. (4,1.7)
    .. controls (2.867,0.567) and (2.3,-0.667) .. (2.3,-2)
    -- (2,-2)
    .. controls (2,-0.867) and (1.333,-0.3) .. (0,-0.3)
    -- cycle;

  \node at (3.3,-1.1) {$\mathcal A_t$};
\end{scope}

\begin{scope}[shift={(7,0)}]
  \draw[thick] (0,0) -- (-2,0);
  \draw[thick] (0,0) -- (0,-2);
  \draw[thick] (0,0) -- (1.6,1.6);

  \node at (-0.85,-0.85) {$0$};
  \node at (-0.85,0.85)  {$x_2$};
  \node at (0.85,-0.85)  {$x_1$};
  \node[right] at (1.05,0.2) {$\Pi_{v=0}$};
\end{scope}

\begin{scope}[shift={(11,-1)}]
  \draw[thick] (0,0) -- (2,0) -- (0,2) -- cycle;

  \node[below left]  at (0,0) {$(0,0)$};
  \node[below]       at (2,0) {$(1,0)$};
  \node[above]       at (0,2) {$(0,1)$};
  \node              at (0.55,0.55) {$\Delta$};
\end{scope}

\end{tikzpicture}
\caption{Left: the rescaled amoeba $\mathcal A_t$ of
$1+z_1+z_2=0$, shown schematically with its three unbounded ends
truncated. Middle: the tropical line of $\max\{0,x_1,x_2\}$;
each chamber is labelled by its dominant term.
Right: the Newton polygon
$\Delta=\operatorname{conv}\{(0,0),(1,0),(0,1)\}$.}
\label{fig:linear-amoeba}
\end{figure}
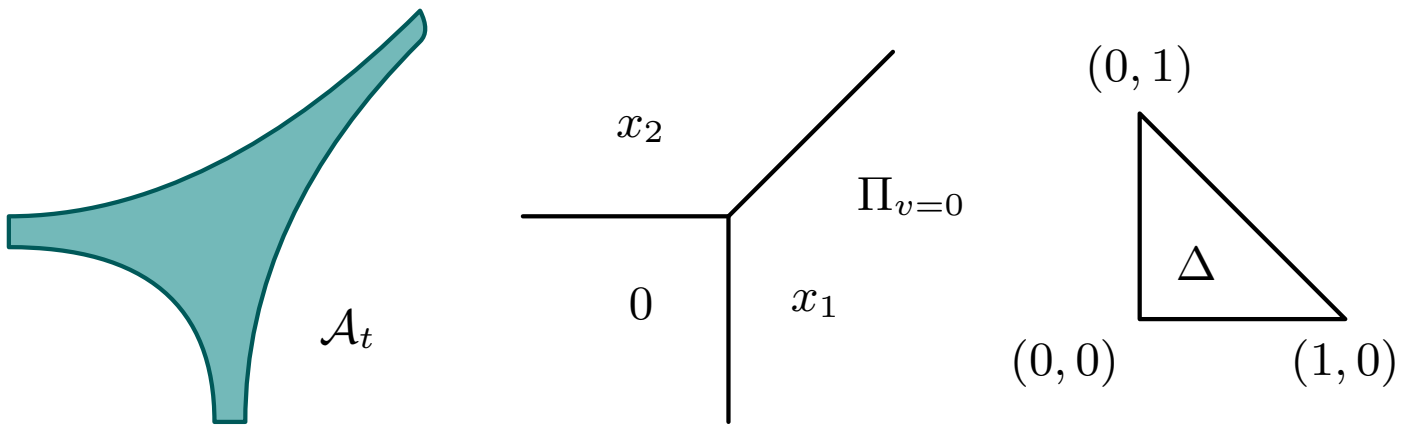

 Define
\[
  L_v(x)=\max_{k\in A}\ell_k(x).
\]
The \emph{tropical hypersurface} \(\Pi_v\) is the locus where this
maximum is attained by at least two indices. The \emph{rescaled amoeba} is defined as the image 
\(\Amo_t=\Log_t(X_t\cap T),\) and Corollary 6.4 in \cite{Mik} asserts   that $$ \mathcal{A}_t \rightarrow \Pi_v,  \  \  \  {\rm when} \ \ t\rightarrow\infty, $$ in the Hausdorff sense. Therefore, $\Pi_v$ is a good approximation of $\mathcal{A}_t$ when $t\gg 1$.

The complement of $\Pi_v$ consists
of the nonempty open convex sets
\[
  C_k=\{x\in\R^n:\ell_k(x)>\ell_j(x)\text{ for all }j\ne k\}.
\]
We call these sets \emph{chambers}. An exponent in the support need
not dominate anywhere, so \(C_k\) may be empty. We also use the
closed polyhedron
\[
  \Xi_k=\{x\in\R^n:\ell_k(x)\ge\ell_j(x)\text{ for all }j\in A\}.
\]

For \(m\in A^\circ\) let
\[
  \Theta_{m,t}=\frac{a_mt^{-v(m)}z^m}{f_t}\,\Omega,
  \qquad
  \eta_{m,t}=\Res_{X_t}\Theta_{m,t}.
\]
To keep the estimates readable, write
\[
  \lambda_{m,t}(p)
   =\CL_\Omega((X_t,\eta_{m,t}),(p,1)),
   \qquad p\in T\setminus X_t.
\]

\begin{lemma}[The boundary of an interior chamber]\label{lem:chamber}
If \(m\in A^\circ\) and \(C_m\ne\varnothing\), then
\(\Xi_m=\overline{C_m}\) is a full-dimensional compact convex
polytope with interior \(C_m\). Its boundary
\(\Gamma_m=\partial\Xi_m\) is contained in \(\Pi_v\) and is
homeomorphic to \(S^{n-1}\).
\end{lemma}

\begin{proof}
Choose \(x_0\in C_m\). If \(y\in\Xi_m\), then every point on the
segment from \(y\) to \(x_0\), other than possibly \(y\), satisfies
all inequalities against \(j\ne m\) strictly. Hence
\(\Xi_m=\overline{C_m}\). Equality against some \(j\ne m\) places
a point on a supporting hyperplane, so \(\Int(\Xi_m)=C_m\).
Such an equality also gives a tie for the maximum, proving
\(\partial\Xi_m\subset\Pi_v\).

Suppose that \(\Xi_m\) were unbounded. Choose \(x_r\in\Xi_m\)
with \(|x_r|\to\infty\) and, after passing to a subsequence,
\(x_r/|x_r|\to u\) with \(|u|=1\). Divide the inequalities for
\(x_r\) by \(|x_r|\) and take limits. We obtain
\[
  \langle m,u\rangle\ge\langle k,u\rangle
  \qquad\text{for all }k\in A.
\]
The linear functional \(y\mapsto\langle y,u\rangle\) would then
attain its maximum on \(\Delta=\Conv(A)\) at the interior point
\(m\), which is impossible. Thus \(\Xi_m\) is bounded. Radial
projection from \(x_0\) identifies the boundary of this
full-dimensional compact convex polytope with \(S^{n-1}\).
\end{proof}

\subsection{The limit theorem}

\begin{theorem}\label{thm:main}
Fix \(m\in A^\circ\). For every nonempty compact set
\(K\subset\R^n\setminus\Pi_v\), there are constants \(B,\delta>0\)
such that, for all sufficiently large \(t\), every \(p\in T\) with
\(\Log_t(p)\in K\) lies outside \(X_t\) and satisfies
\begin{equation}\label{eq:uniform}
  \bigl|\lambda_{m,t}(p)-\ind_{C_m}(\Log_t(p))\bigr|
  \le Bt^{-\delta}.
\end{equation}

If \(C_m\ne\varnothing\), orient
\(\Gamma_m=\partial\Xi_m\) as a boundary. For every family
\(p_t\in T\) such that \(\Log_t(p_t)\to q\notin\Pi_v\),
\begin{equation}\label{eq:link-limit}
  \lim_{t\to\infty}\lambda_{m,t}(p_t)=L(\Gamma_m,q).
\end{equation}
\end{theorem}

\begin{proof}
At each \(x\in K\) there is a unique maximizing index \(j=j(x)\).
The gap between the largest and the second largest of the finitely
many values \(\ell_k(x)\) is continuous and positive on \(K\).
Let \(\delta>0\) be its minimum.

Write \(A_k=a_kt^{-v(k)}p^k\), with \(x=\Log_t(p)\in K\).
By \eqref{eq:size}, for \(k\ne j\),
\[
  \left|\frac{A_k}{A_j}\right|
  \le \left|\frac{a_k}{a_j}\right|t^{-\delta}.
\]
In particular,
\[
  f_t(p)=A_j(1+E),\qquad
  |E|\le D t^{-\delta},\qquad
  D=\max_{j\in A}\sum_{k\ne j}\left|\frac{a_k}{a_j}\right|.
\]
Choose \(t\) large enough that \(Dt^{-\delta}\le1/2\). Then
\(|1+E|\ge1/2\), so \(p\notin X_t\). Notice that this estimate
does not use the arguments of the coordinates of \(p\): no choice
of phases can cancel a term larger than the sum of all the others.

By Proposition~\ref{prop:lift}, \(\lambda_{m,t}(p)=A_m/f_t(p)\).
If \(j=m\), then
\[
  |\lambda_{m,t}(p)-1|
  =\left|\frac{-E}{1+E}\right|\le 2Dt^{-\delta}.
\]
If \(j\ne m\), then
\[
  |\lambda_{m,t}(p)|
  =\left|\frac{A_m/A_j}{1+E}\right|\le 2Dt^{-\delta}.
\]
Since \(j=m\) exactly when \(x\in C_m\), this proves
\eqref{eq:uniform} with \(B=2D\).

Finally, choose a closed ball centred at \(q\) and contained in
its chamber. For all sufficiently large \(t\), this ball contains
\(\Log_t(p_t)\). The chamber indicator is constant on the ball,
so \eqref{eq:uniform} gives the limit \(\ind_{C_m}(q)\).
By Lemma~\ref{lem:chamber} and \eqref{eq:classical}, this is
\(L(\Gamma_m,q)\).
\end{proof}

\begin{remark}
The proof shows directly that
\(\Amo_t\) eventually avoids every compact subset of
\(\R^n\setminus\Pi_v\). This is the only amoeba estimate needed
here. In particular, we do not deduce or assume a two-sided
Hausdorff convergence statement. If \(C_m=\varnothing\),
\eqref{eq:uniform} still holds with limiting value zero, but
there is no full-dimensional \(m\)-chamber whose boundary supplies
the linking interpretation.
\end{remark}

\begin{corollary}[Linear combinations]\label{cor:linear}
Let \(b_m\in\C\) be fixed for \(m\in A^\circ\), and put
\(\eta_t=\sum\limits_{m\in A^\circ}b_m\eta_{m,t}\).
Suppose \(\Log_t(p_t)\to q\in C_j\). Then
\[
  \lim_{t\to\infty}\CL_\Omega((X_t,\eta_t),(p_t,1))
  =
  \begin{cases}
    b_j,&j\in A^\circ,\\
    0,&j\notin A^\circ.
  \end{cases}
\]
\end{corollary}

\begin{proof}
The unique lift of \(\eta_t\) is
\(\sum\limits_{m\in A^\circ}b_m\Theta_{m,t}\). Apply
Theorem~\ref{thm:main} to this finite sum.
\end{proof}

Thus the integer-valued interpretation is a property of the
normalised individual residues. For a linear combination, the
limit instead reads off the coefficient assigned to the chamber.

At a point of \(\Pi_v\) the leading terms may interfere, and
the logarithmic limit need not determine the linking number.
For \(q\in\Pi_v\), put
\[
  I(q)=\{k\in A:l_k(q)=L_v(q)\},\qquad
  F_q(w)=\sum_{k\in I(q)} a_k w^k.
\]

\begin{proposition}\label{prop:wall}
Suppose \(w_t\to w\in T\) and \(F_q(w)\ne0\). Set\\ 
\(p_t=(t^{q_1}w_{t,1},\ldots,t^{q_n}w_{t,n})\).
Then \(p_t\notin X_t\) for all sufficiently large \(t\), and
\[
  \lim_{t\to\infty}
  \CL_\Omega((X_t,\eta_{m,t}),(p_t,1))
  =
  \begin{cases}
    a_m w^m/F_q(w), & m\in I(q),\\
    0,              & m\notin I(q).
  \end{cases}
\]
\end{proposition}

\begin{proof}
Divide the numerator and denominator of \eqref{eq:quotient}
by \(t^{L_v(q)}\). The rescaled denominator tends to \(F_q(w)\)
because a term indexed by \(k\notin I(q)\) has the factor
\(t^{l_k(q)-L_v(q)}\to0\). The same calculation for the
numerator gives the asserted limit.
\end{proof}

\section{A plane cubic: inside, outside, and on the wall}
\label{sec:cubic}

Consider
\[
  f_t(z_1,z_2)=1+z_1^3+z_2^3+t z_1z_2.
\]
Its Newton polygon is
\(\Delta=\Conv\{(0,0),(3,0),(0,3)\}\), and its support consists
of these three vertices and the unique interior lattice point
\(m=(1,1)\). In \(\mathbb{C}P^2\) the equation is
\[
  Z_0^3+Z_1^3+Z_2^3+tZ_0Z_1Z_2=0.
\]
A common zero of the three partial derivatives cannot have a zero
coordinate, since the other two coordinates would then vanish as
well. Multiplying the three derivative equations at a point with
nonzero coordinates gives \(t^3=-27\). Hence \(X_t\) is smooth
for every real \(t>1\).

\subsection{The triangle}

The tropical polynomial and the holomorphic  linking number  are
\[
  L_v(x)=\max\{0,3x_1,3x_2,1+x_1+x_2\},
  \qquad
  \lambda_{m,t}(z)=\frac{t z_1z_2}{f_t(z)}.
\]
Comparing \(1+x_1+x_2\) with each of the other three affine
functions gives
\[
  \Xi_m=
  \left\{x:
  x_1+x_2\ge-1,\quad
  2x_1-x_2\le1,\quad
  2x_2-x_1\le1\right\}.
\]
\begin{figure}[htbp]
\centering
\begingroup
\definecolor{chamberblue}{RGB}{38,91,135}

\begin{minipage}[t]{0.44\linewidth}
\centering
\resizebox{\linewidth}{!}{%
\begin{tikzpicture}[
    x=1.25cm,y=1.25cm,
    font=\small,
    line cap=round,line join=round
]
  \coordinate (O) at (0,0);
  \coordinate (U) at (3,0);
  \coordinate (V) at (0,3);
  \coordinate (M) at (1,1);

  \fill[black!4] (O)--(U)--(V)--cycle;
  \draw[thick] (O)--(U)--(V)--cycle;
  \draw[semithick] (O)--(M) (U)--(M) (V)--(M);

  \fill (O) circle (1.7pt) node[below left]  {$(0,0)$};
  \fill (U) circle (1.7pt) node[below]       {$(3,0)$};
  \fill (V) circle (1.7pt) node[above]       {$(0,3)$};

  \fill[chamberblue] (M) circle (2.6pt);
  \node[above right,text=chamberblue,inner sep=2pt]
    at (M) {$m=(1,1)$};
\end{tikzpicture}%
}
\par\smallskip
{\small\textup{(a)} Newton subdivision}
\end{minipage}
\hfill
\begin{minipage}[t]{0.52\linewidth}
\centering
\resizebox{\linewidth}{!}{%
\begin{tikzpicture}[
    x=1.65cm,y=1.65cm,
    font=\small,
    line cap=round,line join=round,
    ray/.style={semithick,-{Stealth[length=4pt]}},
    boundary/.style={draw=chamberblue,line width=1.4pt},
    oriented/.style={
      postaction={decorate},
      decoration={markings,
        mark=at position 0.65 with
          {\arrow{Stealth[length=5pt,width=5pt]}}}
    }
]
  \coordinate (A) at (-1,0);
  \coordinate (B) at (0,-1);
  \coordinate (C) at (1,1);

  \fill[chamberblue!15] (A)--(B)--(C)--cycle;

  \draw[ray] (A)--(-1.85,0);
  \draw[ray] (B)--(0,-1.65);
  \draw[ray] (C)--(1.60,1.60);

  \draw[boundary] (C)--(A)--(B);
  \draw[boundary,oriented] (B)--(C);

  \fill (0,0) circle (1.7pt)
    node[above,inner sep=2pt] {$q_{\rm in}$};

  \fill (1,0) circle (1.7pt)
    node[right,inner sep=3pt] {$q_{\rm out}$};

  \draw[black,fill=white,line width=0.7pt]
    (-0.5,-0.5) circle (2.1pt)
    node[below left,inner sep=3pt] {$q_{\rm wall}$};
\end{tikzpicture}%
}
\par\smallskip
{\small\textup{(b)} Tropical cubic}
\end{minipage}

\caption{The cubic \(1+z_1^3+z_2^3+t z_1z_2=0\).
Left: the interior exponent \(m=(1,1)\) is incident to three
triangles of the Newton subdivision.
Right: its dual chamber \(C_m\) is the interior of the shaded
triangle. Only the bold, counterclockwise boundary
\(\Gamma_m=\partial\Xi_m\), not the three unbounded rays, enters
the linking number. The marked points are \(q_{\rm in}=(0,0)\),
\(q_{\rm out}=(1,0)\), and
\(q_{\rm wall}=(-\tfrac12,-\tfrac12)\).}
\label{fig:cubic-duality}
\endgroup
\end{figure}
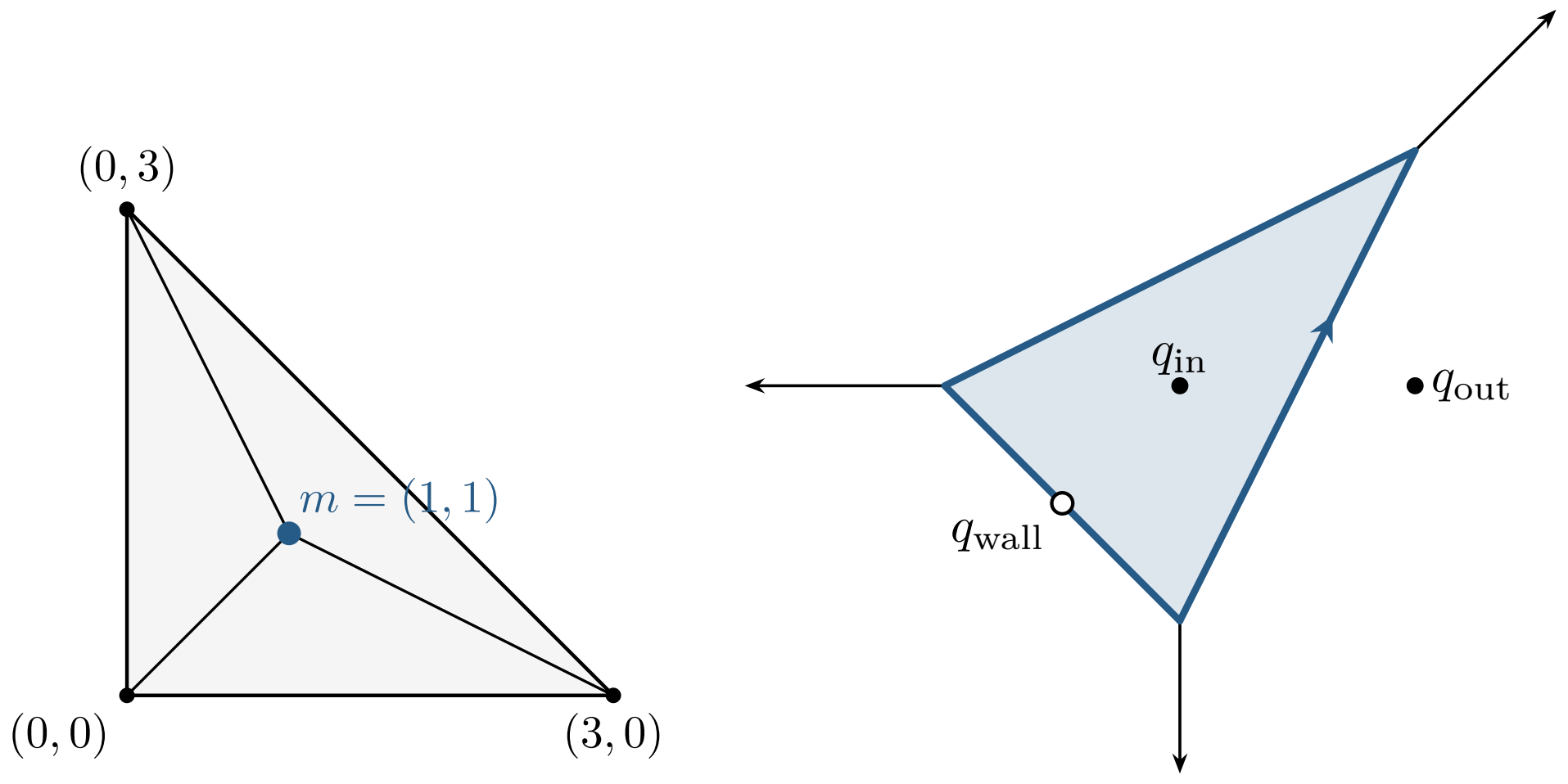
Its vertices, in counterclockwise order, are
\[
  (-1,0),\qquad (0,-1),\qquad (1,1).
\]
The tropical curve consists of the three edges of this triangle
and three unbounded rays: a leftward ray from \((-1,0)\), a
downward ray from \((0,-1)\), and a ray in direction \((1,1)\)
from \((1,1)\). Only the bounded triangle is used for linking.

For \(p_t=(1,1)\), the logarithmic point is \((0,0)\), inside the
triangle, and
\[
  \lambda_{m,t}(p_t)=\frac{t}{t+3}\longrightarrow1.
\]
For \(p_t=(t,1)\), the logarithmic point is \((1,0)\), outside
the triangle, and
\[
  \lambda_{m,t}(p_t)=\frac{t^2}{t^3+t^2+2}\longrightarrow0.
\]

\subsection{A slice through two edges}

One can see both transitions in a single exact calculation.
Take \(p_t=(t^u,t^{-u})\), so that \(\Log_t(p_t)=(u,-u)\).
The line \(x_2=-x_1\) crosses the interior of the triangle for
\(-1/3<u<1/3\). Along the whole line,
\[
  \lambda_{m,t}(p_t)=\frac{t}{1+t^{3u}+t^{-3u}+t}.
\]
Table~\ref{tab:slice} follows by comparing the four terms.
At an edge, the two largest terms have equal magnitude and, for
these positive real coordinates, the same phase.

\begin{table}[htbp]
\centering
\caption{The linking limit along the line \(x_2=-x_1\), with
positive real point coordinates. The two edge values depend on
this choice of phases.}
\label{tab:slice}
\small
\begin{tabular}{@{}lll@{}}
\toprule
Range of \(u\) & Largest term(s) in the denominator
              & Limit of \(\lambda_{m,t}\)\\
\midrule
\(u<-1/3\)       & \(t^{-3u}\)          & \(0\)\\
\(u=-1/3\)       & \(t^{-3u}\) and \(t\)& \(1/2\)\\
\(-1/3<u<1/3\)   & \(t\)                & \(1\)\\
\(u=1/3\)        & \(t^{3u}\) and \(t\) & \(1/2\)\\
\(u>1/3\)        & \(t^{3u}\)           & \(0\)\\
\bottomrule
\end{tabular}
\end{table}

\subsection{The same wall point, different limits}

Is \(1/2\) a natural boundary value of the linking function?
Fix the midpoint \(q=(-1/2,-1/2)\) of the third edge. For any
\(s\in\C\) with \(|s|=1\) and \(s\ne-1\), set
\[
  p_t=(t^{-1/2},s\,t^{-1/2}).
\]
These points all have \(\Log_t(p_t)=q\), not merely
\(\Log_t(p_t)\to q\). Nevertheless,
\[
  t z_1z_2=s,\qquad
  f_t(p_t)=1+s+(1+s^3)t^{-3/2},
\]
and therefore
\[
  \lim_{t\to\infty}\lambda_{m,t}(p_t)=\frac{s}{1+s}.
\]
For \(s=1\) the limit is \(1/2\); for \(s=i\) it is
\((1+i)/2\). At the excluded phase \(s=-1\), even the lower
order terms cancel and \(p_t\in X_t\) for every \(t\), so the
linking function is not defined there.

Off the tropical curve, the logarithmic point determines the
limit regardless of phase. On this edge, the logarithmic point is
fixed and the phase alone changes the answer. This is why the
classical linking interpretation stops at the chamber boundary.

\section*{Acknowledgement of AI assistance}

The author used ChatGPT 6.2 (OpenAI) to assist with the organisation
and drafting of the manuscript and with revisions to proofs and
examples. The author takes responsibility for the mathematical
claims, references, and final text.

\end{document}